\documentclass{article}

\usepackage{url}

\usepackage{amsmath, amssymb, amstext, amsopn, amsxtra, graphicx, color, calligra, hyperref, url, ulem, polynom, permute, ifthen, float, epstopdf, linsys, marvosym, rotating}

\begin{document}
\setlength{\parskip}{2mm}

\newcommand{\field}[1]{\ensuremath{\mathbb #1}}
\newcommand{\MM}{\field M}
\renewcommand{\vec}[1]{\ensuremath{\begin{pmatrix}#1\end{pmatrix}}}
\newcommand{\ur}[1]{\ensuremath{\mathrm{ur}(#1)}}
\newcommand{\pcomment}[1]{\mbox{#1}\hskip 15mm}
\newcommand{\imply}{\ensuremath{\Longrightarrow} }
\newcommand{\nin}{\not\in}
\newcommand{\inv}{^{-1}}
\newcommand{\qar}{\begin{eqnarray*}}
\newcommand{\raq}{\end{eqnarray*}}
\newcommand{\qed}{\hskip 1cm {\bf QED}}
\newcommand{\pf}{\em{Proof}  }

\renewcommand{\em}{\textit}

%\title{UNIT REGULAR ELEMENTS IN CORNER RINGS}
\title{Unit regular elements in corner rings}

\author{T.\ Y.\ LAM\(^*\) \hskip 1in WILL MURRAY\\ \\
\textit {University of California, Berkeley, CA  94720, U.S.A.}}
\date{\textit{(Received on 10 August, 1996; Accepted on 3 November, 1996)}}

\maketitle

\begin{abstract}
\noindent For any ring \(R\), some characterizations are obtained for unit regular elements in a corner ring \(eRe\) in terms of unit regular elements in \(R\).\\


\end{abstract}

\section{Introduction}

Throughout this note, \(R\) denotes a ring with identity.  An element \(a \in R\) is said to be \em{unit regular} in \(R\) if \(a = aua\) for some \(u \in U(R)\), the group of units of \(R\).  Let us write \ur R\ for the set of all unit regular elements in \(R\).  If \(\ur R = R\), \(R\) is said to be a \em{unit regular ring} (see \cite{ehrlich1}).

Unit regular rings play an important role in the structure theory of (von Neumann) regular rings since they are precisely those regular rings over which the finitely generated (say left) projective modules satisfy the cancellation law [3: (4.5)].  Alternatively, unit regular rings can also be characterized as regular rings of stable range one [3: (4.12)].

For any idempotent \(e \in R\), the ring \(eRe\) (with identity \(e\)) is said to be a \em{corner ring} of \(R\).  It is easy to show that, for any \(a \in eRe\):
\[
\pcomment{(*)} a \in \ur{eRe} \imply a \in \ur R.
\]
(Take \(b = 0\) in the implication (1) \imply (4) in the Theorem below.)  However, the \em{converse} of (*) is not true in general.  To construct a counterexample, let \(S\) be a ring with an element \(s\) which is regular but not unit regular.  Choose \(t \in S\) such that \(s = sts\), and let \(R = \MM_2(S), e = \vec{1 & 0\\ 0 & 0} \in R\).  Then, for \(a = \vec{s & 0\\ 0 & 0} \in eRe\) and \(u = \vec{t & 1\\ 1 & 0} \in R\), we have 
\[
aua = \vec{s & 0\\ 0 & 0}\vec{t & 1\\ 1 & 0}\vec{s & 0\\ 0 & 0} = \vec{sts & 0\\ 0 & 0} = a.
\]
Since \(u \in U(R)\) (with inverse \vec{0 & 1\\ 1 & -t}), we have \(a \in \ur R\).  However \(a \nin \ur{eRe}\) since \(eRe \cong S\) and \(s \nin \ur S\).

The fact that (*) is only a one-way implication leaves one to wonder whether there exist a characterization of the elements of  \ur{eRe} in terms of those of \ur R.  In this note, we'll provide some such characterizations.  As an easy consequence of these, for instance, we get a direct ``elementwise'' proof for the fact that \em{any corner ring of a unit regular ring is unit regular}.  Some of the implications in the main result below can also be proved by using the technique of endomorphism rings in \cite{ehrlich2}.  However, we believe a direct proof without recourse to the results of \cite{ehrlich2} is of independent interest.  Another advantage of our proof is that it applies equally well to 1-sided unit regular elements and 1-sided unit regular rings.  

We thank K.\ R.\ Goodearl for pointing out to us that the converse of (*) is false in general.  This remark of Goodearl prompted us to find an ``iff'' statement to replace the ``only if'' statement (*).  

\section{Main results}

We come now to the main result of this paper.

\noindent{\sc Theorem}   Let \(e = e^2 \in R, f = 1-e\), and \(a \in eRe\).  Then the following statements are equivalent:

\noindent (1)  \(a \in \ur{eRe}\);\\
(2)  \(a + f \in \ur R\);\\
(3)  \(a + b \in \ur R\) for any \(b \in U(fRf)\);\\
(3')  \(a + b \in \ur R\) for some \(b \in U(fRf)\);\\
(4)  \(a + b \in \ur R\) for any \(b \in \ur{fRf}\);\\
(5)  \(a + b \in \ur R\) for some \(b \in fRf\) which is neither a left nor a right 0-divisor of \(fRf\).

\pf  Clearly, (4) \imply (3) \imply (2) \imply (3') \imply (5), so it suffices for us to prove (1) \imply (4) and (5) \imply (1).

\noindent (1) \imply (4).  Let \(D = eRe \times fRf\), which is a subring of \(R\) (with the same identity \(e + f = 1\)).  For any \(b \in \ur{fRf}\), we have clearly
\[
a+b \in \ur{eRe} \times \ur{fRf} = \ur D \subseteq \ur R.
\]
\noindent (5) \imply (1).  For \(b\) as in (5), write \(a+b = (a+b)u(a+b)\) where \(u \in U(R)\).  Then
\[
a+b = aua + bub + bua + aub,
\]
where the four summands on the RHS lie in \(eRe, fRf, fRe\), and \(eRf\) respectively.  By the Peirce Decomposition Theorem [5: p.318], we have \(a = aua\), \(b=bub\) and \(bua = aub = 0\).  Now let \(v = u\inv \in R\), and consider the two elements
\[
v' := eve \in eRe, \hskip 5mm u' := e(u-ubu)e \in eRe.
\]
We have
\qar
au'a & = & ae(u-ubu)ea\\
& = & a(u-ubu)a\\
& = & aua - (aub)ua\\
& = & aua\\
& = & a,
\raq
so it suffices to show that \(u', v'\) are mutually inverse in the ring \(eRe\).  Now \((1-bu)f \in fRf\) and \([(1-bu)f]b = b-bub = 0\), so \((1-bu)f = 0\) (since \(b\) is not a right 0-divisor in \(fRf\)), and hence \((1-bu)e = 1-bu\).  Using this, we get
\qar
u'v' & = & eu(1-bu)e \cdot eve\\
& = & eu(1-bu)ve\\
& = & e(1-ub)uve\\
& = & e(1-ub)e\\
& = & e,
\raq
since \(be \in fRf \cdot e = 0\).  Similarly, the fact that \(b\) is not a left 0-divisor in \(fRf\) implies that \(e(1-ub) = 1-ub\), and using this, we can calculate as before that \(v'u' = e\).  Therefore, \(u' \in U(eRe)\), and the fact that \(a = au'a\) shows that \(a \in \ur{eRe}\), as desired.  \qed

\noindent{\sc Remarks}

\noindent (A)  There is an obvious statement (4') obtained by replacing ``for any'' in (4) by ``for some''.  However, (4') is in general \em{weaker} than the other conditions, as we have seen before (by taking \(b=0\)!).

\noindent (B)  The proof of (5) \imply (1) actually shows more; namely, if \(a+b\) is a \em{right unit regular} element of \(R\) (in the sense of Ehrlich \cite{ehrlich2}) for some \(b \in fRf\) which is not a right 0-divisor of \(fRf\), then \(a\) is a \em{right unit regular} element of the corner ring \(eRe\).

\noindent(C)  It is also of interest to point out that the proof of (5) \imply (1) did not really make full use of the fact that \(uv = 1 = vu\).  In fact, as the reader may check, all we needed in the proof was that \(uve = e = evu\).  Therefore, (1) would have followed if we just assume that, for the element \(u\) above, there exists \(v \in R\) such that \(uv - 1 \in Rf\) and \(vu - 1 \in fR\).

As an immediate consequence of the Theorem above, we obtain:  

\noindent{\sc Corollary}  Let \(e = e^2 \in R\).  If \(R\) is a unit regular ring, so is \(eRe\).

\pf  Let \(f = 1-e\).  For any \(a \in eRe\), we have \(a + f \in \ur R\).  By (2) \imply (1) in the above Theorem, we have \(a \in \ur{eRe}\).  (A similar proof, using Remark (B) above, shows that if \(R\) is a 1-sided unit regular ring in the sense of Ehrlich \cite{ehrlich2}, then so is \(eRe\) for any idempotent \(e \in R\).)  \qed

The Corollary above is a well-known result in the theory of unit regular rings.  However, the only known proof of it depends on a characterization theorem of Ehrlich \cite{ehrlich2} and Handelman \cite{handelman} (in terms of module cancellation) for an endomorphism ring to be unit regular.  (See also [6:  Ex.\ (21.9)].)  Our treatment above provided the first ``elementwise'' proof of this corner ring result.

\noindent \(^*\)The author to whom correspondence should be addressed.  [fax: (510)642-8204]

\vfill\eject
\end{document}